\documentclass{article}

\title{The maximum modulus of a trigonometric trinomial}
\author{Stefan Neuwirth}
\date{}

\usepackage[british]{babel}

\usepackage{amssymb}
\usepackage{amsmath}
\usepackage{mathrsfs}
\usepackage{amsthm}
\usepackage{graphicx}

\theoremstyle{plain}
\newtheorem{thm}{Theorem}[section]
\newtheorem*{Ch}{N.~G. Chebotar\"ev's formula}
\newtheorem{lem}[thm]{Lemma}
\newtheorem{prp}[thm]{Proposition}
\newtheorem{cor}[thm]{Corollary}

\theoremstyle{definition}

\newtheorem*{Not}{Notation}
\newtheorem*{acknowledgements}{Acknowledgement}
\newtheorem{ep}[thm]{Extremal problem}

\theoremstyle{remark}
\newtheorem{rem}[thm]{Remark}

\renewcommand{\le}{\leqslant}
\renewcommand{\ge}{\geqslant}

\def\Cont{\mathscr C}

\def\th{\vartheta}

\def\rh{\varrho}

\def\la{\lambda}

\def\eps{\epsilon}

\def\ie{\textit{i.\,e.}}

\def\e{\mkern1mu\mathrm e\mkern1mu}
\def\iu{\mkern1mu\mathrm i\mkern1mu}
\def\ei#1{\e^{\iu #1}} 
\def\emi#1{\e^{-\iu #1}}

\providecommand{\abs}[1]{\lvert#1\rvert}
\providecommand{\bigabs}[1]{\bigl\lvert#1\bigr\rvert}

\def\io[#1,#2]{\mathopen]#1,\allowbreak#2\mathclose[}
\def\iog[#1,#2]{\mathopen]#1,\allowbreak#2]}
\def\bigiog[#1,#2]{\bigl]#1,\allowbreak#2\bigr]}
\def\bigiod[#1,#2]{\bigl[#1,\allowbreak#2\bigr[}
\def\iod[#1,#2]{[#1,\allowbreak#2\mathclose[}

\def\Z{\mathbb Z}

\def\R{\mathbb R}

\begin{document}

\maketitle

\begin{abstract}

  Let $\Lambda$ be a set of three integers and let $\Cont_\Lambda$ be the
  space of $2\pi$-periodic functions with spectrum in $\Lambda$
  endowed with the maximum modulus norm. We isolate the maximum modulus points $x$ of
  trigonometric trinomials $T\in\Cont_\Lambda$ and prove that $x$ is unique
  unless $\abs T$ has an axis of symmetry. This permits to compute the exposed
  and the extreme points of the unit ball of $\Cont_\Lambda$, to describe how
  the maximum modulus of $T$ varies with respect to the arguments of its Fourier coefficients
  and to compute the norm of unimodular relative Fourier multipliers on
  $\Cont_\Lambda$. We obtain in particular the Sidon constant of $\Lambda$.
\end{abstract}

\section{Introduction}

Let $\la_1$, $\la_2$ and $\la_3$ be three pairwise distinct integers. Let
$r_1$, $r_2$ and $r_3$ be three positive real numbers. Given three real
numbers $t_1$, $t_2$ and $t_3$, let us consider the \emph{trigonometric trinomial}
\begin{equation}
  \label{tt}
  T(x)=r_1\ei{(t_1+\la_1x)}+r_2\ei{(t_2+\la_2x)}+r_3\ei{(t_3+\la_3x)}
\end{equation}
for $x\in\R$. The $\la$'s are the \emph{frequencies} of the trigonometric
trinomial $T$, the $r$'s are the \emph{moduli} or \emph{intensities} and the
$t$'s the \emph{arguments} or \emph{phases} of its \emph{Fourier
  coefficients} $r_1\ei{t_1}$, $r_2\ei{t_2}$ and $r_3\ei{t_3}$.

\begin{figure}[htbm]
  \centering
  \includegraphics{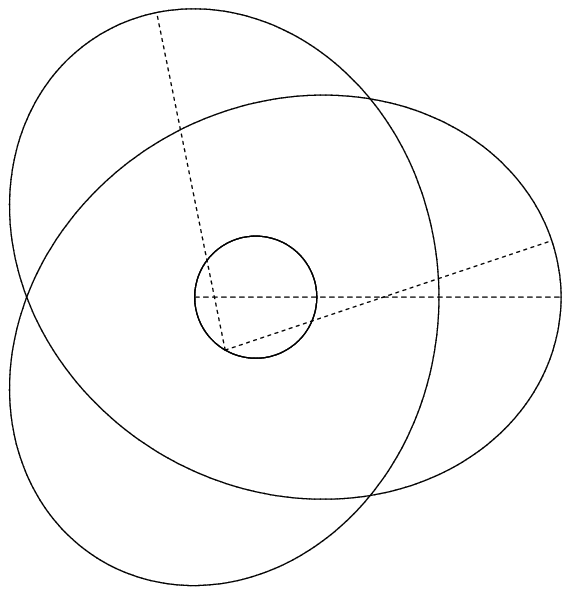}
  \caption{The unit circle, the hypotrochoid $z=4\emi{2x}+\ei{x}$, its unique point at maximum distance to $-1$ and its two points at maximum distance to $-\ei{\pi/3}$.}
  \label{fig:1}
\end{figure}

The maximum modulus of a trigonometric trinomial has a geometric
interpretation. Without loss of generality, we may assume that $\la_2$ is
between $\la_1$ and $\la_3$.  Let $H$ be the curve with complex equation
\begin{equation}
  \label{hypo}
z=r_1\ei{(t_1-(\la_2-\la_1)x)}+r_3\ei{(t_3+(\la_3-\la_2)x)}\quad
(-\pi<x\le\pi)\text.
\end{equation}
$H$ is a \emph{hypotrochoid}: it is drawn by a point at distance $r_3$ to the
centre of a circle with radius $r_1\abs{\la_2-\la_1}/\abs{\la_3-\la_2}$ that
rolls inside another circle with radius
$r_1\abs{\la_3-\la_1}/\abs{\la_3-\la_2}$. The maximum modulus of \eqref{tt} is
the maximum distance of points $z\in H$ to a given point $-r_2\ei{t_2}$ of the
complex plane.  Figure~\ref{fig:1} illustrates a particular case.

We deduce an interval on which $T$ attains its maximum modulus
independently of the moduli of its Fourier coefficients (see
Theorem~\ref{mmpttt}\,$(a)$ for a detailed answer.) We prove in
particular the following result.

\begin{thm}
  \label{unique}
  Let $d=\gcd(\la_2-\la_1,\allowbreak\la_3-\la_2)$ and let
  $\tau$ be the distance of
  \begin{equation}\label{tau}
    \frac{\la_2-\la_3}d\allowbreak t_1
    +\frac{\la_3-\la_1}d\allowbreak t_2
    +\frac{\la_1-\la_2}d\allowbreak t_3
  \end{equation}
  to $2\pi\Z$. The trigonometric trinomial $T$ attains its maximum modulus at a
  unique point modulo $2\pi/d$, with multiplicity $2$, unless $\tau=\pi$. 
\end{thm}

Theorem~\ref{unique} shows that if there are two points of the hypotrochoid $H$ at maximum distance to $-r_2\ei{t_2}$, it is so only because
$-r_2\ei{t_2}$ lies on an axis of symmetry of $H$.

We obtain a precise description of those trigonometric trinomials that attain
their maximum modulus twice modulo $2\pi/d$ (see Theorem~\ref{mmpttt}\,$(c)$.)
Their r\^ole is illustrated by the following geometric result. Let us first
put up the proper functional analytic framework. Let
$\Lambda=\{\la_1,\la_2,\la_3\}$ be the \emph{spectrum} of the trigonometric
trinomial $T$ and denote $\e_\la\colon x\mapsto\ei{\la x}$.  Let
$\Cont_\Lambda$ be the space of functions spanned by the $\e_{\la}$ with
$\la\in\Lambda$, endowed with the maximum modulus norm. Recall that a point $P$ of a
compact convex set $K$ is \emph{exposed} if there is a hyperplane that meets
$K$ only in $P$; $P$ is \emph{extreme} if it is not the midpoint of any two
other points of $K$. 
\begin{thm}
  \label{extrpts}
  Let $K$ be the unit ball of the space $\Cont_\Lambda$ and let $P\in
  K$.
  \begin{enumerate}
  \item The point $P$ is an exposed point of $K$ if and only if $P$ is either a trigonometric monomial
    $\ei\alpha\e_\la$ with $\alpha\in\R$ and $\la\in\Lambda$ or a
    trigonometric trinomial that attains its maximum modulus, $1$, at
    two points modulo $2\pi/d$. Every linear functional on $\Cont_\Lambda$
    attains its norm on an exposed point of $K$.
  \item The point $P$ is an extreme point of $K$ if and only if $P$ is either
    a trigonometric monomial $\ei\alpha\e_\la$ with $\alpha\in\R$ and
    $\la\in\Lambda$ or a trigonometric trinomial such that $1-{\abs P}^2$ has
    four zeroes modulo $2\pi/d$, counted with multiplicities.
  \end{enumerate}
\end{thm}

We describe the dependence of the maximum modulus of the trigonometric
trinomial $T$ on the arguments.  The general issue has been addressed for a
long time. The articles \cite{li24,sa33} are two early references.

Let us state our main theorem, that solves an elementary case of the complex
Man\-del$'$\-shtam problem, a term coined in \cite{dm90}. It appeared
originally in electrical circuit theory: ``L.~I. Mandel$'$shtam communicated
me a problem on the phase choice of electric currents with different
frequencies such that the capacity of the resulting current to blow is
minimal'' \cite[p.~396]{ce49}.

\begin{ep}[Complex Man\-del$'$\-shtam problem]
  To find the minimum of the maximum modulus of a trigonometric polynomial with
  given Fourier coefficient moduli. 
\end{ep}

\begin{thm}
  \label{max-decreasing}
  The maximum modulus of~\eqref{tt} 
  is a strictly decreasing function of $\tau$. In
  particular, 
  \begin{multline*}
    \min_{t_1,t_2,t_3}\max_x
    \bigabs{r_1\ei{(t_1+\la_1x)}+r_2\ei{(t_2+\la_2x)}+r_3\ei{(t_3+\la_3x)}}\\
    =\max_{x} 
    \bigabs{\epsilon_1r_1\ei{\la_1x} + \epsilon_2r_2\ei{\la_2x} + \epsilon_3r_3\ei{\la_3x}}
  \end{multline*}
  if $\epsilon_1$, $\epsilon_2$ and $\epsilon_3$ are real signs $+1$ or
  $-1$ such that $\epsilon_i\epsilon_j=-1$, where $i,j,k$ is a
  permutation of $1,2,3$ such that the power of $2$ in $\la_i-\la_j$ is
  greater than the power of $2$ in $\la_i-\la_k$ and in
  $\la_k-\la_j$.
\end{thm}

Our result shows that the maximum modulus is minimal when the phases are
chosen in opposition, independently of the intensities $r_1$, $r_2$ and $r_3$.

The decrease of the maximum modulus of~\eqref{tt} may be bounded as shown in
the next result.

\begin{thm}
  \label{how-dec}
  Suppose that $\la_2$ is between $\la_1$ and $\la_3$.  The quotient of the
  maximum modulus of~\eqref{tt} by $\bigabs{r_1+r_2\ei{\tau
      d/\abs{\la_3-\la_1}}+r_3}$ is a strictly increasing function of $\tau$
  unless $r_1:r_3=\abs{\la_3-\la_2}:\abs{\la_2-\la_1}$, in which case it is
  constantly equal to $1$.
\end{thm}

\begin{figure}[htbm]
  \centering
  \includegraphics{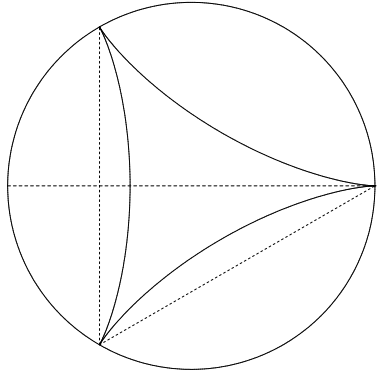}
  \caption{The unit circle, the deltoid $z=(1/3)\emi{2x}+(2/3)\ei{x}$, its
    unique point at maximum distance to $-1$ and its two points at maximum
    distance to $-\ei{\pi/3}$.}
  \label{fig2}
\end{figure}

When $r_1:r_3=\abs{\la_3-\la_2}:\abs{\la_2-\la_1}$, the hypotrochoid $H$ with
equation~\eqref{hypo} is a \emph{hypocycloid} with $\abs{\la_3-\la_1}/d$
cusps: the rolling point is \emph{on} the rolling circle. Figure~\ref{fig2}
illustrates a particular case.

We may deduce from Theorem~\ref{how-dec} a less precise but handier
inequality.

\begin{thm}
  \label{th-4}
  Let 
  \begin{equation}
    \label{D}
    D=\frac
    {\max(\abs{\la_2-\la_1},\abs{\la_3-\la_2},\abs{\la_3-\la_1})}
    {\gcd(\la_2-\la_1,\la_3-\la_1)}
  \end{equation}
  be the quotient of the diameter of $\Lambda$ by $d$.
Let $t'_1$, $t'_2$ and $t'_3$ be another three real numbers and define
  similarly $\tau'$.  If $\tau>\tau'$, then
  \begin{multline*}
    \max_x\bigabs{r_1\ei{(t_1+\la_1x)}+r_2\ei{(t_2+\la_2x)}+r_3\ei{(t_3+\la_3x)}}\\
    \ge\frac{\cos(\tau /2D)}{\cos(\tau'/2D)}
    \max_x\bigabs{r_1\ei{(t'_1+\la_1x)}+r_2\ei{(t'_2+\la_2x)}+r_3\ei{(t'_3+\la_3x)}}
  \end{multline*}
  with equality if and only if
  $r_1:r_2:r_3=\abs{\la_3-\la_2}:\abs{\la_3-\la_1}:\abs{\la_2-\la_1}$.
\end{thm}

If we choose $\tau'=0$ in the last result, we get the solution to an
elementary case of the following extremal problem.
\begin{ep}
  To find the minimum of the maximum modulus of a trigonometric polynomial with
  given spectrum, Fourier coefficient arguments and moduli sum.
\end{ep}

\begin{thm}
  \label{th-3}
  We have
  \begin{equation*}
    \frac
    {\max_x\bigabs{r_1\ei{(t_1+\la_1x)}+r_2\ei{(t_2+\la_2x)}+r_3\ei{(t_3+\la_3x)}}}
    {r_1+r_2+r_3}
    \ge\cos(\tau /2D)
  \end{equation*}
  with equality if and only if $\tau=0$ or
  $r_1:r_2:r_3=\abs{\la_3-\la_2}:\abs{\la_3-\la_1}:\abs{\la_2-\la_1}$.
\end{thm}

The dependence of the maximum modulus of~\eqref{tt} on the arguments may also
be expressed as properties of relative multipliers.  Given three real numbers
$t_1$, $t_2$ and $t_3$, the linear operator on $\Cont_\Lambda$ defined by
$\e_{\la_j}\mapsto\ei{t_j}\e_{\la_j}$ is a \emph{unimodular
  relative Fourier multiplier}: it multiplies each Fourier coefficient of
elements of $\Cont_\Lambda$ by a fixed unimodular number; let us denote it by
$(t_1,t_2,t_3)$. Consult \cite{ha87} for general background on relative
multipliers.

\begin{thm}
  \label{multiplier}
  The unimodular relative
  Fourier multiplier $(t_1,t_2,t_3)$ has norm
  \begin{equation*}
    \cos\bigl((\pi-\tau)/2D\bigr)\big/\cos(\pi/\allowbreak2D),
  \end{equation*}
  where $\tau$ is defined as in Theorem~\ref{unique} and $D$ is given by~\eqref{D}, and attains its norm exactly
  at functions of the form
  \begin{equation*}
    r_1\ei{(u_1+\la_1x)}+r_2\ei{(u_2+\la_2x)}+r_3\ei{(u_3+\la_3x)}
  \end{equation*}
  with $r_1:r_2:r_3=\abs{\la_3-\la_2}:\abs{\la_3-\la_1}:\abs{\la_2-\la_1}$ and
  \begin{equation*}
    \frac{\la_2-\la_3}du_1+\frac{\la_3-\la_1}du_2+\frac{\la_1-\la_2}du_3=\pi\mod{2\pi}.
  \end{equation*}
\end{thm}

The maximum of the norm of unimodular relative Fourier multipliers is the
\emph{complex unconditional constant} of the canonical basis
$(\e_{\la_1},\e_{\la_2},\e_{\la_3})$ of $\Cont_\Lambda$. As
\begin{equation*}
  r_1+r_2+r_3=\max_x\bigabs{r_1\ei{\la_1x}+r_2\ei{\la_2x}+r_3\ei{\la_3x}},
\end{equation*}
this constant is the minimal constant $C$ such that
\begin{equation*}
  r_1+r_2+r_3\le C\max_x\bigabs{r_1\ei{(u_1+\la_1x)}+r_2\ei{(u_2+\la_2x)}+r_3\ei{(u_3+\la_3x)}};
\end{equation*}
it is therefore the \emph{Sidon constant} of $\Lambda$. It is also the
solution to the following extremal problem.
\begin{ep}[Sidon constant problem]
  To find the minimum of the maximum modulus of a trigonometric polynomial
  with given spectrum and Fourier coefficient moduli sum.
\end{ep}

We obtain the following result, setting $\tau=\pi$ in Theorem~\ref{multiplier}.

\begin{cor}
  The Sidon constant of $\Lambda$ is $\sec(\pi d/\allowbreak2D)$. It is
  attained exactly at functions of the form given in Theorem~\ref{multiplier}.
\end{cor}

Let us now give a brief description of this article. In Sections~\ref{sec2} and~\ref{sec2bis}, we use carefully the invariance of the maximum modulus under rotation,
translation and conjugation to reduce the arguments $t_1$, $t_2$ and $t_3$ of
the Fourier coefficients of the trigonometric trinomial $T$ to the variable
$\tau$.  Section~\ref{sec3} shows how to further reduce this study to the
trigonometric trinomial
\begin{equation}
  \label{Tred}
  r_1\emi{kx}+r_2\ei{t}+r_3\ei{lx}
\end{equation}
with $k$ and $l$ positive coprime integers and $t\in[0,\pi/(k+l)]$. In
Section~\ref{secloc}, we prove that~\eqref{Tred} attains its maximum modulus
for $x\in[-t/k,t/l]$. Section~\ref{sec4} studies the variations of the modulus
of~\eqref{Tred} for $x\in[-t/k,t/l]$: it turns out that it attains its
absolute maximum only once on that interval. This yields Theorem~\ref{unique}.
Section~\ref{secmm} restates the results of the two previous sections for a
general trigonometric trinomial $T$. Section~\ref{sec:exp} is dedicated to the
proof of Theorem~\ref{extrpts}. In Section~\ref{seccheb}, we compute the
directional derivative of the maximum modulus of~\eqref{Tred} with respect to
the argument $t$ and prove Theorems~\ref{max-decreasing}, \ref{how-dec},
\ref{th-4} and~\ref{th-3}. In Section~\ref{secnm}, we prove
Theorem~\ref{multiplier} and show how to lift unimodular relative Fourier
multipliers to operators of convolution with a linear combination of two Dirac
measures.  Section~\ref{secsidon} replaces our computation of the Sidon
constant in a general context; it describes the initial motivation for this
research.

Part of these results appeared previously, with a different proof, in
\cite[Cha\-pter II.10]{ne99} and in \cite{ne01b}.

\begin{Not}
  Throughout this article, $\la_1$, $\la_2$ and $\la_3$ are three pairwise
  distinct integers, $\Lambda=\{\la_1,\la_2,\la_3\}$ and
  $d=\gcd(\la_2-\la_1,\allowbreak\la_3-\la_2)$. If $\la$ is an integer,
  $\e_\la$ is the function $x\mapsto\ei{\la x}$ of the real variable $x$.  A
  \emph{trigonometric polynomial} is a linear combination of functions $\e_\la$; it is a
  \emph{monomial}, \emph{binomial} or \emph{trinomial} if this linear
  combination has one, two or three nonzero coefficients, respectively. The
  normed space $\Cont_\Lambda$ is the three-dimensional space of complex
  functions spanned by $\e_\la$ with $\la\in\Lambda$, endowed with the maximum
  modulus norm. The Dirac measure $\delta_x$ is the linear functional
  $T\mapsto T(x)$ of evaluation at $x$ on the space of continuous functions.
\end{Not}

\section{Isometric relative Fourier multipliers}
\label{sec2}

The r\^ole of Quantity~\eqref{tau} is explained by the following lemma.

\begin{lem}
  \label{isom}
  Let $t_1$, $t_2$ and $t_3$ be real numbers. The unimodular relative Fourier
  multiplier $M=(t_1,\allowbreak t_2,\allowbreak t_3)$ is an isometry on
  $\Cont_\Lambda$ if and only if
  \begin{equation}
    \label{th}
    \frac{\la_2-\la_3}dt_1
    +\frac{\la_3-\la_1}dt_2
    +\frac{\la_1-\la_2}dt_3
    \in2\pi\Z:
  \end{equation}
  it is a unimodular multiple of a translation: there are real numbers
  $\alpha$ and $v$ such that $Mf(x)=\ei\alpha f(x-v)$ for all
  $f\in\Cont_\Lambda$ and all $x\in\R$.
\end{lem}

\begin{proof}
  If $M$ is a unimodular multiple of a
  translation by a real number $v$, then
  \begin{equation*}
    \bigabs{r_1\ei{(t_1+\la_1v)}+r_2\ei{(t_2+\la_2v)}+r_3\ei{(t_3+\la_3v)}}=r_1+r_2+r_3,
  \end{equation*}
  which holds if and only if
  \begin{equation}
    \label{thx}
    t_1+\la_1v=t_2+\la_2v=t_3+\la_3v\quad\text{modulo $2\pi$.}
  \end{equation}
  There is a $v$ satisfying~\eqref{thx} if and only if
  Equation~\eqref{th} holds as~\eqref{thx} means that there exist
  integers $a_1$ and $a_3$ such that
  \begin{equation*}
    v
    =\frac{t_2-t_1+2\pi a_1}{\la_1-\la_2}
    =\frac{t_2-t_3+2\pi a_3}{\la_3-\la_2}.
  \end{equation*}
  If $t_1$, $t_2$ and $t_3$ are three real numbers satisfying~\eqref{th}, let
  $v$ be such that~\eqref{thx} holds. Then
  \begin{multline*}
    {r_1\ei{(t_1+u_1+\la_1 x)}
      +r_2\ei{(t_2+u_2+\la_2 x)}
      +r_3\ei{(t_3+u_3+\la_3 x)}}\\
    =\ei{(t_2+\la_2v)}\left({r_1\ei{(u_1+\la_1(x-v))}+r_2\ei{(u_2+\la_2(x-v))}+r_3\ei{(u_3+\la_3(x-v))}}\right)
  \end{multline*}
  for all real numbers $u_1$, $u_2$, $u_3$ and $x$.
\end{proof}

\section{The arguments of the Fourier coefficients of a trigonometric trinomial}
\label{sec2bis}

A translation and a rotation permit to reduce the three arguments of the
Fourier coefficients of a trigonometric trinomial to just one variable. Use of the involution
$\overline{f(-x)}$ of $\Cont_\Lambda$ permits to restrain even further
the domain of that variable.

\begin{lem}
  \label{upto}
  Let $t_1$, $t_2$ and $t_3$ be real numbers and let $\tilde t_2$ be the representant of
  \begin{equation}
    \label{tilde}
    \frac{\la_2-\la_3}{\la_3-\la_1}\allowbreak t_1
    +t_2
    +\frac{\la_1-\la_2}{\la_3-\la_1}\allowbreak t_3 
  \end{equation}
  modulo $2\pi/\abs{\la_3-\la_1}$ in $\bigiod[-\pi d/\abs{\la_3-\la_1},\pi d/\abs{\la_3-\la_1}]$.
  \begin{enumerate}
  \item There are real numbers $\alpha$ and $v$ such that
    \begin{multline}
      \label{nett}
      {r_1\ei{(t_1+\la_1x)}+r_2\ei{(t_2+\la_2x)}+r_3\ei{(t_3+\la_3x)}}\\
      =\ei\alpha\left({r_1\ei{\la_1(x-v)}+r_2\ei{(\tilde
          t_2+\la_2(x-v))}+r_3\ei{\la_3(x-v))}}\right)
    \end{multline}
    for all $x$.
  \item Let $t=\abs{\tilde t_2}$ be the distance of~\eqref{tilde} to $(2\pi
    d/\abs{\la_3-\la_1})\Z$. There is a sign $\eps\in\{+1,-1\}$ such that
    \begin{multline*}
      \bigabs{r_1\ei{(t_1+\la_1x)}+r_2\ei{(t_2+\la_2x)}+r_3\ei{(t_3+\la_3x)}}\\
      =\bigabs{r_1\ei{\la_1\eps(x-v)}+r_2\ei{(t+\la_2\eps(x-v))}+r_3\ei{\la_3\eps(x-v))}}
    \end{multline*}
    for all $x$.
  \end{enumerate}
\end{lem}

\begin{proof}
  $(a)$. The argument $\tilde t_2$ is chosen so that the relative multiplier
  $(t_1,t_2-\tilde t_2,t_3)$ is an isometry. 

  $(b)$. If $\tilde t_2$ is negative, take the conjugate under the modulus of
  the right hand side in~\eqref{nett}.
\end{proof}

\begin{rem}
  We have also
  \begin{multline*}
    \bigabs{r_1\ei{\la_1x} +r_2\ei{(t+2\pi d/(\la_3-\la_1)+\la_2x)}
      +r_3\ei{\la_3x}}\\
    =\bigabs{r_1\ei{\la_1(x-v)} +r_2\ei{(t+\la_2(x-v))}
      +r_3\ei{\la_3(x-v)}}
  \end{multline*}
  for all $x$ and $t$, where $v$ satisfies $\la_1v=2\pi
  d/(\la_3-\la_1)+\la_2v=\la_3v$ modulo $2\pi$, that is
  \begin{equation*}
    v=\frac{2m\pi}{\la_3-\la_1}\quad\text{with $m$ an inverse of $\frac{\la_3-\la_2}d$
      modulo $\frac{\la_3-\la_1}d$.}
  \end{equation*}
\end{rem}

\section{The frequencies of a trigonometric trinomial}\label{sec3}

We may suppose without loss of generality that $\la_1<\la_2<\la_3$.
Let $k=(\la_2-\la_1)/\allowbreak d$ and $l=(\la_3-\la_2)/d$. Then
\begin{equation*}
  \bigabs{r_1\ei{\la_1x}
    +r_2\ei{(t+\la_2x)}
    +r_3\ei{\la_3x}}
  =\bigabs{r_1\emi{k(dx)}
    +r_2\ei{t}
    +r_3\ei{l(dx)}}.
\end{equation*}
A homothety by $d^{-1}$  permits to restrain our study to the
function
\begin{equation*}
  f(t,x)
  ={\bigabs{r_1\emi{kx}
      +r_2\ei{t}
      +r_3\ei{lx}}}^2
\end{equation*}
for $x\in\R$ with $k$ and $l$ two positive coprime numbers and
$t\in\bigl[0,\pi/(k+l)\bigr]$.  We have
\begin{align}
  {f(-t,x)}
  &={f(t,-x)}
  \label{pair}\\
  {f\bigl(t+2\pi/(k+l),x\bigr)}
  &={f(t,x-2m\pi/(k+l))}
  \label{per}
\end{align}
for all $x$ and $t$, where 
$m$ is the inverse of $l$ modulo $k+l$.
In particular, if $t=\pi/(k+l)$, we have the symmetry relation
\begin{equation}
  \label{symmetry}
  {f\bigl(\pi/(k+l),x\bigr)}
  ={f\bigl(\pi/(k+l),2m\pi/(k+l)-x\bigr)}\text.
\end{equation}

\section{Location of the maximum point}
\label{secloc}

The purpose of our first proposition is to deduce a small interval on
which a trigonometric trinomial attains its maximum modulus. Note that
a trigonometric binomial attains its maximum modulus at a point that
depends only on the phase of its coefficients:
\begin{itemize}
\item $\bigabs{r_1\emi{kx}+r_2\ei{t}}$ attains its maximum at $-t/k$
  independently of $r_1$ and $r_2$,
\item $\bigabs{r_1\emi{kx}+r_3\ei{lx}}$ attains its maximum at $0$
  independently of $r_1$ and $r_3$,
\item $\bigabs{r_2\ei{t}+r_3\ei{lx}}$ attains its maximum at $t/l$
  independently of $r_2$ and $r_3$.
\end{itemize}
The next proposition shows that if
the point at which a trigonometric trinomial attains its maximum
modulus changes with the intensity of its coefficients, it changes
very little; we get bounds for this point that are independent of
the intensities.

\begin{prp}
  \label{loc}
  Let $k,l$ be two positive coprime integers. Let $r_1$, $r_2$ and
  $r_3$ be three positive real numbers. Let
  $t\in\bigl[0,\pi/(k+l)\bigr]$.  Let
  \begin{equation*}
    f(x)={\bigabs{r_1\emi{kx}+r_2\ei{t}+r_3\ei{lx}}}^2
  \end{equation*}
  for $x\in\R$.
  \begin{enumerate}
  \item The function $f$ attains its absolute maximum in the
    interval $[-t/k,\allowbreak t/l]$.
  \item If $ f$ attains its absolute maximum at a point $y$ outside of
    $[-t/k,t/l]$ modulo $2\pi$, then $t=\pi/(k+l)$ and
    $2m\pi/(k+l)-y$ lies in $[-t/k,t/l]$ modulo $2\pi$, where $m$ is
    the inverse of $l$ modulo $k+l$.
  \end{enumerate}
\end{prp}

\begin{proof}
  $(a)$. We have
  \begin{multline}
    \label{f(x)carre}
    {f(x)}=r_1^2+r_2^2+r_3^2+2\cdot\\
    \bigl(r_1r_2\cos(t+kx)+r_1r_3\cos\bigl((k+l)x\bigr)+r_2r_3\cos(t-lx)\bigr).
  \end{multline}
  Let us prove that ${f}$ attains its absolute maximum on
  $[-t/k,\allowbreak t/l]$.  Let $y$ be outside of $[-t/k,t/l]$ modulo
  $2\pi$. Let $I$ be the set of all $x\in[-t/k,t/l]$ such that
  \begin{equation*}
    \left\{
      \begin{aligned}
        \cos(t+kx)&\ge\cos(t+ky)\\
        \cos\bigl((k+l)x\bigr)&\ge\cos\bigl((k+l)y\bigr)\\
        \cos(t-lx)&\ge\cos(t-ly).
      \end{aligned}
    \right.
  \end{equation*}
  Note that if $x\in[-t/k,t/l]$, then
  \begin{equation*}
    \left\{
      \begin{aligned}
        t+kx&\in[0,(k+l)t/l]\\ 
        (k+l)x&\in[-(k+l)t/k,(k+l)t/l]\\ 
        t-lx&\in[0,(k+l)t/k],
      \end{aligned}
    \right.
  \end{equation*}
  and that $(k+l)t/k$, $(k+l)t/l\subset[0,\pi]$. Let
  \begin{itemize}
  \item $\alpha$ be the distance of $t/k+y$ to $(2\pi/k)\Z$,
  \item $\beta$ be the distance of $y$ to $\bigl(2\pi/(k+l)\bigr)\Z$,
  \item $\gamma$ be the distance of $t/l-y$ to $(2\pi/l)\Z$.
  \end{itemize}
  Then
  \begin{equation}
    \label{II}
    I=[-t/k,t/l]\cap
    [-t/k-\alpha,-t/k+\alpha]\cap
    [-\beta,\beta]\cap
    [t/l-\gamma,t/l+\gamma].
  \end{equation}
  Let us check that $I$ is the nonempty interval
  \begin{equation}\label{I}
    I=\bigl[\max(-t/k,-\beta,t/l-\gamma),\min(t/l,-t/k+\alpha,\beta)\bigr].
  \end{equation}
  In fact, we have the following triangular inequalities:
  \begin{itemize}
  \item $-\beta\le-t/k+\alpha$ because $t/k$ is the
    distance of $(t/k+y)-y$ to $(2\pi/k(k+l))\Z$;
  \item $t/l-\gamma\le-t/k+\alpha$ because $t/l+t/k$ is the distance of
    $(t/k+y)+(t/l-y)$ to $(2\pi/kl)\Z$;
  \item $t/l-\gamma\le\beta$ because
    $t/l$ is the distance of
    $(t/l-y)+y$ to $(2\pi/l(k+l))\Z$.
  \end{itemize}
  The other six inequalities that are necessary to deduce~\eqref{I}
  from~\eqref{II} are obvious. 

  $(b)$. We have proved in $(a)$ that there is an $x\in[-t/k,t/l]$ such
  that $\cos(t+kx)\ge\cos(t+ky)$,
  $\cos\bigl((k+l)x\bigr)\ge\cos\bigl((k+l)y\bigr)$ and
  $\cos(t-lx)\ge\cos(t-ly)$. In fact, at least one of these inequalities is
  strict unless there are signs $\delta,\eps,\eta\in\{-1,1\}$ such that
  $t+kx=\delta(t+ky)$, $t-lx=\eps(t-ly)$ and $(k+l)x=\eta(k+l)y$
  modulo $2\pi$. Two out of these three signs are equal and the corresponding
  two equations imply the third one with the same sign. This system is
  therefore equivalent to
  \begin{equation*}
    \left\{
      \begin{aligned}
        k(x-y)&=0\\
        l(x-y)&=0
      \end{aligned}
    \right.
    \quad\textrm{or}\quad
    \left\{
      \begin{aligned}
        k(x+y)&=-2t\\
        l(x+y)&=2t
      \end{aligned}
    \right.
  \end{equation*}
  modulo $2\pi$. The first pair of equations yields $x=y$ modulo
  $2\pi$ because $k$ and $l$ are coprime. Let $m$ be an inverse of $l$
  modulo $k+l$; then the second pair of equations is equivalent to
  \begin{equation*}
    \left\{
      \begin{aligned}
        2(k+l)t&=0\\
        x+y&=2mt
      \end{aligned}
    \right.
  \end{equation*}
  modulo $2\pi$. Therefore $g$ does not attain its absolute maximum at
  $y$ unless $t=\pi/(k+l)$ and $2m\pi/(k+l)-y\in[-t/k,t/l]$.
\end{proof}

\begin{rem}
  This proposition is a complex version of \cite[Lemma~2.1]{re95}.  
\end{rem}

\section{Uniqueness of the maximum point}
\label{sec4}

Note that
\begin{align*}
  r_1\emi{kx}+r_2\ei{t}+r_3\ei{lx}
  &=r_3\emi{l(-x)}+r_2\ei{t}+r_1\ei{k(-x)}\\
  &=r'_1\emi{k'x'}+r_2\ei{t}+r'_3\ei{l'x'}
\end{align*}
with $r'_1=r_3$, $r'_3=r_1$, $k'=l$, $l'=k$ and $x'=-x$.
We may therefore suppose without loss of generality that $kr_1\le lr_3$.

Our second proposition studies the points at which a trigonometric trinomial
attains its maximum modulus. Note that if $k=l=1$, the derivative of
${\abs{f}}^2$ has at most $4$ zeroes, so that the modulus of $f$ has at most
two maxima and attains its absolute
maximum in at most two points. Proposition~\ref{prp-2} shows that this is true
in general, and that if it may attain its absolute maximum in two points, it
is so only because of the symmetry given by~\eqref{symmetry}.

\begin{prp}
  \label{prp-2}
  Let $k,l$ be two positive coprime integers. Let $r_1$, $r_2$ and
  $r_3$ be three positive real numbers such that $kr_1\le lr_3$. Let
  $t\in\iog[0,\pi/(k+l)]$. Let
  \begin{equation*}
    f(x)={\bigabs{r_1\emi{kx}+r_2\ei{t}+r_3\ei{lx}}}^2
  \end{equation*}
  for $x\in[-t/k,t/l]$.
  \begin{enumerate}
  \item There is a point $x^*\in[0,t/l]$ such that $\mathrm df/\mathrm dx>0$
    on $\io[-t/k,x^*]$ and $\mathrm df/\mathrm dx<0$ on $\io[x^*,t/l]$.
  \item There are three cases:
    \begin{enumerate}
    \item ${f}$ attains its absolute maximum at $0$ if and only if
      $kr_1=lr_3$;
    \item ${f}$ attains its absolute maximum at $t/l$ if and only if
      $l=1$, $t=\pi/(k+1)$ and $k^2r_1r_2+(k+1)^2r_1r_3-r_2r_3\le0$;
    \item otherwise, ${f}$ attains its absolute maximum in
      $\io[0,t/l]$.
    \end{enumerate}
  \item The function $f$ attains its absolute maximum with
    multiplicity $2$ unless $l=1$, $t=\pi/(k+1)$ and
    $k^2r_1r_2+(k+1)^2r_1r_3-r_2r_3=0$, in which case it attains its
    absolute maximum at $\pi/(k+1)$ with multiplicity $4$.
  \end{enumerate}
\end{prp}
\begin{proof}
  $(a)$. By Proposition~\ref{loc}, the derivative of $f$ has a zero in
  $[-t/k,t/l]$. Let us study the sign of the derivative of ${f}$. Equation~\eqref{f(x)carre} yields
  \begin{equation}\label{df(x)}
    \frac12\dfrac{\mathrm d{f}}{\mathrm dx}(x)
    =-kr_1r_2\sin(t+kx)-(k+l)r_1r_3\sin\bigl((k+l)x\bigr)+lr_2r_3\sin(t-lx).
  \end{equation}

  We wish to compare $\sin(t+kx)$ with $\sin(t-lx)$: note that
  \begin{equation*}
    \sin(t+kx)-\sin(t-lx)
    =2\sin\bigl((k+l)x/2\bigr)
    \cos\bigl(t+(k-l)x/2\bigr).
  \end{equation*}
  Suppose that $x\in[-t/k,t/l]$. Then
  \begin{gather*}
    -\pi\le-\pi/k\le-(k+l)t/k\le(k+l)x\le(k+l)t/l\le\pi/l\le\pi\\
    0\le2t+(k-l)x\le
    \begin{cases}
      2t+(l-k)t/k=(k+l)t/k&\text{if $k\le l$}\\
      2t+(k-l)t/l=(k+l)t/l&\text{if $l\le k$}
    \end{cases}
    \le\pi.
  \end{gather*}

  Suppose that $x\in\iod[-t/k,0]$: then it follows that
  $\sin(t+kx)\le\sin(t-lx)$ and $\sin\bigl((k+l)x\bigr)\le0$ with
  equality if and only if $k=1$ and $-x=t=\pi/(1+l)$. This yields
  with $kr_1\le lr_3$ that
  \begin{equation}
    \label{1}
    \frac12\dfrac{\mathrm d{f}}{\mathrm dx}(x)
    \ge-(k+l)r_1r_3\sin\bigl((k+l)x\bigr)\ge0
  \end{equation}
  with equality if and only if $k=1$ and $-x=t=\pi/(1+l)$. 

  If $l\ge2$ and $x\in[0,t/l]$, then
  \begin{equation*}
    \left\{
      \begin{aligned}
        t+kx&\in[t,(k+l)t/l]\subset[t,\pi/2]\\ 
        (k+l)x&\in[0,(k+l)t/l]\subset[t,\pi/2]\\ 
        t-lx&\in[0,t]\subset[0,\pi/3],
      \end{aligned}
    \right.
  \end{equation*}
  so that the second derivative of ${f}$ is strictly negative
  on $[0,t/l]$: its derivative is strictly decreasing on
  this interval and $(a)$ is proved.
  Suppose that $l=1$ and consider $g(x)=f(t-x)$ for $x\in[0,t]$: we have
  to prove that there is a point $x^*$ such that $\mathrm dg/\mathrm dx>0$ on
  $\io[0,x^*]$ and $\mathrm dg/\mathrm dx<0$ on $\io[x^*,t]$. Put
  $\alpha=(k+1)t$: then
  \begin{equation*}
    \frac12\dfrac{\mathrm d{g}}{\mathrm dx}(x)
    =kr_1r_2\sin(\alpha-kx)
    +(k+1)r_1r_3\sin\bigl(\alpha-(k+1)x\bigr)
    -r_2r_3\sin x
  \end{equation*}
  and it suffices to prove that
  \begin{equation}
    \label{g}
    \dfrac{1}{2\sin x}\dfrac{\mathrm d{g}}{\mathrm dx}(x)
    =kr_1r_2\frac{\sin(\alpha-kx)}{\sin x}
    +(k+1)r_1r_3\frac{\sin\bigl(\alpha-(k+1)x\bigr)}{\sin x}
    -r_2r_3
  \end{equation}
  decreases strictly with $x\in\iog[0,{\alpha/(k+1)}]$.  Let us study
  the sign of
  \begin{equation*}
    \frac{\mathrm d}{\mathrm dx}\frac{\sin(\alpha-kx)}{\allowbreak\sin x}
    =\frac{-k\cos(\alpha-kx)\sin x-\sin(\alpha-kx)\cos x}{\sin^2x}
  \end{equation*}
  for $\alpha\in\iog[0,\pi]$ and $x\in\iog[0,{\alpha/k}]$. If $k=1$,
  then
  \begin{equation*}
    -k\cos(\alpha-kx)\sin x-\sin(\alpha-kx)\cos x=-\sin\alpha\le0
  \end{equation*}
  and the inequality is strict unless $\alpha=\pi$. Let us prove by
  induction on $k$ that $k\cos(\alpha-kx)\sin x+\sin(\alpha-kx)\cos
  x>0$ for all $k\ge2$, $\alpha\in\iog[0,\pi]$ and
  $x\in\iog[0,{\alpha/k}]$. This will end the proof of $(a)$. Let
  $k\ge1$ and $x\in\iog[0,{\alpha/(k+1)}]$. Then
  \begin{multline*}
    (k+1)\cos\bigl(\alpha-(k+1)x\bigr)\sin x
    +\sin\bigl(\alpha-(k+1)x\bigr)\cos x\\
    \begin{aligned}
     &=\begin{aligned}[t]
        (k+1)\cos(\alpha-kx)\cos x\sin x+(k+1)\sin(\alpha-kx)\sin^2 x\\
        +\sin(\alpha-kx)\cos^2 x-\cos(\alpha-kx)\sin x\cos x&
      \end{aligned}\\
      &=\begin{aligned}[t]
        \bigl(k\cos(\alpha-kx)\sin x+\sin(\alpha-kx)\cos x\bigr)\cos x\\
        +(k+1)\sin(\alpha-kx)\sin^2 x&
      \end{aligned}\\
      &\ge(k+1)\sin(\alpha-kx)\sin^2 x>0
    \end{aligned}
  \end{multline*}

  $(b)$. 1. By Proposition~\ref{loc} and $(a)$, $f$ attains its absolute maximum at $0$ if and only if $0$ is a
  critical point for $f$. We have
  \begin{equation*}
    \frac12\dfrac{\mathrm d{f}}{\mathrm dx}(0)
    =(lr_3-kr_1)r_2\sin t\ge0
  \end{equation*}
  and equality holds if and only if $kr_1=lr_3$. 

  2. We have
  \begin{equation*}
    \frac12\dfrac{\mathrm d{f}}{\mathrm dx}(t/l)
    =\bigl(-kr_1r_2-(k+l)r_1r_3\bigr)\sin\bigl((k+l)t/l\bigr)\le0
  \end{equation*}
  and equality holds if and only if $l=1$ and $t=\pi/(k+1)$.  Let
  $l=1$ and $t=\pi/(k+1)$ and let us use the notation introduced in
  the last part of the proof of $(a)$: we need to characterise the case that
  $g$ has a maximum at $0$. As $\alpha=\pi$ and
  \begin{equation}\label{d2g}
    \frac1{2\sin x}\dfrac{\mathrm d{g}}{\mathrm dx}(x)
    =k^2r_1r_2+(k+1)^2r_1r_3-r_2r_3+o(x)
  \end{equation}
  decreases strictly with $x\in\iog[0,{\pi/(k+1)}]$, $g$ has a maximum at
  $0$ if and only if $k^2r_1r_2+(k+1)^2r_1r_3-r_2r_3\le0$.

  $(c)$. If $l\ge2$, then the second derivative of $f$ is strictly
  negative on $[0,t/l]$. If $l=1$, then the derivative of~\eqref{g} is
  strictly negative on $\iog[0,\alpha/(k+1)]$: this yields that the
  second derivative of $g$ can only vanish at $0$. By $(b)\,2$, $g$ has
  a maximum at $0$ only if $t=\pi/(k+1)$; then
  \begin{gather}
    \frac12\dfrac{\mathrm d^2{g}}{\mathrm dx^2}(0)
    =k^2r_1r_2+(k+1)^2r_1r_3-r_2r_3
    \label{d2g2}\\
    \frac12\dfrac{\mathrm d^4{g}}{\mathrm dx^4}(0)
    =-k^4r_1r_2-(k+1)^4r_1r_3+r_2r_3
    \label{d2g3}
  \end{gather}
  If~\eqref{d2g2} vanishes, then the sum of~\eqref{d2g3} with
  \eqref{d2g2} yields
  \begin{equation*}
    \frac12\dfrac{\mathrm d^4{g}}{\mathrm dx^4}(0)
    =-k(k+1)r_1\bigl((k-1)kr_2+(k+1)(k+2)r_3\bigr)<0.\qedhere
  \end{equation*}
\end{proof}

\begin{rem}
  We were able to prove directly that the system
  \begin{equation*}
    \left\{
      \begin{aligned}
        &{f(x)}={f(y)}\\
        &\dfrac{\mathrm d{f}}{\mathrm dx}(x)   
        =\dfrac{\mathrm d{f}}{\mathrm dx}(y)=0\\
        &\dfrac{\mathrm d^2{f}}{\mathrm dx^2}(x),   
        \dfrac{\mathrm d^2{f}}{\mathrm dx^2}(y)\le0
      \end{aligned}
    \right.
  \end{equation*}
  implies $x=y$ modulo $2\pi$ or $t=\pi/(k+l)$ and
  $x+y=2m\pi/(k+l)$, but our computations are very involved and
  opaque.
\end{rem}

\begin{rem}
  Suppose that $l=k=1$. If $t\in\io[0,\pi/2]$, it is necessary to
  solve a generally irreducible quartic equation in order to compute
  the maximum of $f$. If $t=\pi/2$, it suffices to solve a
  linear equation and one gets the following expression for
  \begin{math}
        \max_{x}{\bigabs{r_1\emi x+\iu r_2+r_3\ei x}}:
  \end{math}
  \begin{equation*}
    \begin{cases}
      (r_1+r_3)\sqrt{1+{r_2^2}/{4r_1r_3}}&\text{if $\bigabs{r_1^{-1}-r_3^{-1}}<4r_2^{-1}$}\\
      r_2+\abs{r_3-r_1}&\text{otherwise.}
    \end{cases}
  \end{equation*}
  In the first case, the maximum is attained at the two points $x^*$ such
  that $\sin x^*=r_2(r_3-r_1)/4r_1r_3$.
\end{rem}

\begin{rem}
  Suppose that $l=1$ and $k=2$. If $t\in\io[0,\pi/3]$, it is
  necessary to solve a generally irreducible sextic equation in order to
  compute the maximum of $f$. If $t=\pi/3$, it suffices to
  solve a quadratic equation and one gets the following expression for
  \begin{math}
    \max_{x}{\bigabs{r_1\emi{2x}+r_2\ei{\pi/3}+r_3\ei x}}:
  \end{math}
  if $r_1^{-1}-4r_3^{-1}<9r_2^{-1}$, then its square makes
  \begin{equation*}
    r_1^2+\frac23r_2^2+r_3^2+r_1r_2+2r_1r_3\biggl[{\biggl({\Bigl(\frac{r_2}{3r_3}\Bigr)}^2+\frac{r_2}{3r_1}+1\biggr)}^{3/2}-{\Bigl(\frac{r_2}{3r_3}\Bigr)}^3\biggr]
  \end{equation*}
  and the maximum is attained at the two points $x^*$ such that
  \begin{equation*}
    2\cos(\pi/3-x^*)={\biggl({\Bigl(\frac{r_2}{3r_3}\Bigr)}^2+\frac{r_2}{3r_1}+1\biggr)}^{1/2}-\frac{r_2}{3r_3};
  \end{equation*}
  otherwise, it makes $-r_1+r_2+r_3$.
\end{rem}

\section{The maximum modulus points of a trigonometric trinomial}
\label{secmm}

If we undo all the reductions that we did in Sections~\ref{sec2bis}
and~\ref{sec3} and in the beginning of Section~\ref{sec4}, we get the following theorem.

\begin{thm}
  \label{mmpttt}
  Let $\la_1$, $\la_2$ and $\la_3$ be three pairwise distinct integers such
  that $\la_2$ is between $\la_1$ and $\la_3$. Let $r_1$, $r_2$ and $r_3$ be
  three positive real numbers.  Given three real numbers $t_1$, $t_2$
  and $t_3$, consider the trigonometric trinomial
  \begin{equation}
    \label{T}
    T(x)={{r_1\ei{(t_1+\la_1x)}+r_2\ei{(t_2+\la_2x)}+r_3\ei{(t_3+\la_3x)}}}
  \end{equation}
  for $x\in\R$. Let $d=\gcd(\la_2-\la_1,\allowbreak\la_3-\la_2)$
  and choose integers $a_1$ and $a_3$ such that
  \begin{equation*}
    \tau=\frac{\la_2-\la_3}d(t_1-2\pi a_1)+\frac{\la_3-\la_1}dt_2+\frac{\la_1-\la_2}d(t_3-2\pi a_3),
  \end{equation*}
  satisfies $\abs\tau\le\pi$. Let $\tilde t_1=t_1-2\pi a_1$ and
  $\tilde t_3=t_3-2\pi a_3$.
  \begin{enumerate}
  \item The trigonometric trinomial $T$ attains its maximum modulus at
    a unique point of the interval bounded by
    \begin{math}
      (\tilde t_1-t_2)/\allowbreak(\la_2-\la_1)
    \end{math}
    and
    \begin{math}
      (t_2-\tilde t_3)/\allowbreak(\la_3-\la_2).
    \end{math}
    More precisely,
    \begin{itemize}
    \item if $r_1\abs{\la_2-\la_1}\le r_3\abs{\la_3-\la_2}$, then this point
      is between
      \begin{math}
        (\tilde t_1-\tilde t_3)/(\la_3-\la_1)
      \end{math}
      and
      \begin{math}
        (t_2-\tilde t_3)/(\la_3-\la_2);
      \end{math}
    \item if $r_1\abs{\la_2-\la_1}\ge r_3\abs{\la_3-\la_2}$, then this point
      is between
      \begin{math}
        (\tilde t_1-\tilde t_3)/(\la_3-\la_1)
      \end{math}
      and
      \begin{math}
        (\tilde t_1-t_2)/(\la_2-\la_1);
      \end{math}
    \item $T$ attains its maximum modulus at
      $(\tilde t_1-\tilde t_3)/(\la_3-\la_1)$ if and only if
      $r_1\abs{\la_2-\la_1}\allowbreak=\allowbreak r_3\abs{\la_3-\la_2}$ or
      $\tau=0$.
    \end{itemize}
  \item The function $T$ attains its maximum modulus at a unique point
    modulo $2\pi/d$, and with multiplicity $2$, unless $\abs\tau=\pi$.
  \item Suppose that $\abs\tau=\pi$, \ie,
    \begin{equation}
      \label{opp}
      \frac{\la_2-\la_3}dt_1+\frac{\la_3-\la_1}dt_2+\frac{\la_1-\la_2}dt_3=\pi\mod{2\pi}\text.
    \end{equation}
    Let $s$ be a solution to $2t_1+\la_1s=2t_2+\la_2s=2t_3+\la_3s$ modulo
    $2\pi$: $s$ is unique modulo $2\pi/d$. 
    Then ${T(s-x)}=\ei{(2t_2+\la_2s)}\overline{T(x)}$ for all $x$. Suppose that
    $\abs{\la_3-\la_2}\le\abs{\la_2-\la_1}$. There are three cases.
    \begin{enumerate}
    \item If $\la_2-\la_1=k(\la_3-\la_2)$ with $k\ge2$ 
      integer and
      \begin{equation*}
        r_1^{-1}-k^2r_3^{-1}\ge(k+1)^2r_2^{-1},
      \end{equation*}
      then ${T}$ attains its maximum modulus, $-r_1+r_2+r_3$, only at
      $x=(t_2-\tilde t_3)/\allowbreak(\la_3-\la_2)$ modulo $2\pi/d$,
      with multiplicity $2$ if the inequality is strict and with
      multiplicity $4$ if there is equality;
    \item if $\la_2-\la_1=\la_3-\la_2$ and
      \begin{equation*}
      \bigabs{r_1^{-1}-r_3^{-1}}\ge4r_2^{-1},
    \end{equation*}
    then ${T}$ attains its maximum modulus, $r_2+\abs{r_3-r_1}$, at a
    unique point $x$ modulo $2\pi/d$, with multiplicity $2$ if the
    inequality is strict and with multiplicity $4$ if there is
    equality. This point is $(t_2-\tilde
    t_3)/\allowbreak(\la_3-\la_2)$ if $r_1<r_3$, and $(\tilde
    t_1-t_2)/\allowbreak(\la_2-\la_1)$ if $r_3<r_1$;
    \item otherwise ${T}$ attains its maximum modulus at exactly two points
      $x$ and $y$ modulo $2\pi/d$, with multiplicity $2$, where $x$ is strictly between
      \begin{math}
        (\tilde t_1-t_2)/\allowbreak(\la_2-\la_1)
      \end{math}
      and
      \begin{math}
        (t_2-\tilde t_3)/\allowbreak(\la_3-\la_2),
      \end{math}
      and $x+y=s$ modulo $2\pi/d$.
    \end{enumerate}
    Note that $s-x=x$ modulo $2\pi/d$ in Cases 1 and 2.
  \end{enumerate}
\end{thm}

\section{Exposed and extreme points of the unit ball of
  $\Cont_\Lambda$}
\label{sec:exp}

The characterisation of the maximum modulus points of a trigonometric
trinomial permits to compute the exposed and the extreme points of the
unit ball of $\Cont_\Lambda$. We begin with a lemma.
\begin{lem}
  \label{deter}
  \begin{enumerate}
  \item A trigonometric trinomial with a given spectrum that attains its
    maximum modulus at two given points modulo $2\pi/d$ is determined by its value at these points.
  \item The trigonometric trinomials with a given spectrum that attain their
    maximum modulus with multiplicity $4$ at a given point and have a given
    value at this point lie on a parabola.
  \end{enumerate}
\end{lem}

\begin{proof}
  We will use the notation of Theorem~\ref{mmpttt}. Without loss of
  generality, we may suppose that $\la_1=-k$, $\la_2=0$ and $\la_3=l$ with $k$
  and $l$ two positive coprime integers. Let $x$ and $y$ be two real numbers
  that are different modulo $2\pi/d$, let $\th$ and $\zeta$ be
  real numbers and let $\rh$ be a positive real number. 

  $(a)$. Let us prove that at most one trigonometric trinomial $T$ attains its
  maximum modulus at $x$ and $y$ and satisfies $T(x)=\rh\ei\th$ and
  $T(y)=\rh\ei\zeta$. Let us translate $T$ by $(x+y)/2$: we may suppose that
  $x+y=0$. Let us divide $T$ by $\ei{(\th+\zeta)/2}$: we may suppose that
  $\th+\zeta=0$. As $T$ attains its maximum modulus at the two points $x$ and
  $y$, we have $s=x+y=0$ and $2t_1-ks=2t_2=2t_3+ls=\th+\zeta=0$ modulo $2\pi$.
  Therefore $t_1=t_2=t_3=0$ modulo $\pi$. Let $p_j=\ei{t_j}r_j$: the $p_j$ are
  nonzero real numbers. We have
  \begin{equation*}
    T(x)=p_1\emi{kx}+p_2+p_3\ei{lx}=\rh\ei{\th}\text,
  \end{equation*}
  so that, multiplying by $\emi{\th}$ and taking real and imaginary parts,
  \begin{align}
    p_1\cos(\th+kx)+p_3\cos(\th-lx)&=\rh-p_2\cos\th
    \label{cy}\\
    p_1\sin(\th+kx)+p_3\sin(\th-lx)&=-p_2\sin\th\text.
    \label{sy1}
  \end{align}
  The computation
  \begin{equation*}
    \frac12\frac{\mathrm d{\abs T}^2}{\mathrm dx}(x)
    =\Re\Bigl(\overline{T(x)}\frac{\mathrm dT}{\mathrm dx}(x)\Bigr)
    =\Re\Bigl(\overline{T(x)}\bigl(-\iu k p_1 \emi{kx} + \iu l p_3 \ei{lx}\bigr)\Bigr),
  \end{equation*}
   yields
  \begin{equation}\label{sy2}
    kp_1\sin(\th+kx)-lp_3\sin(\th-lx)=0.
  \end{equation}
  Equations~\eqref{sy1} and~\eqref{sy2} yield $p_1$ and $p_3$ as linear
  functions of $p_2$ because $\sin(\th+kx)\allowbreak\sin(\th-lx)\ne0$:
  otherwise both factors vanish, so that $\th=x=0$ modulo $\pi$ and $x=y$
  modulo $2\pi$. As $\rh\ne0$, Equation~\eqref{cy} has at most one solution in
  $p_2$.

  $(b)$. We may suppose that $l=1$. Let us determine all trigonometric
  trinomials $T$ that attain their maximum modulus at $x$ with multiplicity
  $4$ and satisfy $T(x)=\rh\ei\th$.  Let us translate $T$ by $x$: we may
  suppose that $x=0$. Let us divide $T$ by $\ei{\th}$: we may suppose that
  $\th=0$. As $T$ attains its maximum modulus at $0$ with multiplicity $4$, we
  have $s-0=0$ and $2t_1-ks=2t_2=2t_3+s=2\th=0$ modulo $2\pi$. Therefore
  $t_1=t_2=t_3=0$ modulo $\pi$. 
  Let $p_j=\ei{t_j}r_j$: the $p_j$ are nonzero real numbers and satisfy the system
  \begin{equation*}
    \left\{
    \begin{aligned}
    &p_1+p_2+p_3=\rh\\
    &k^2p_1p_2+{(k+1)}^2p_1p_3+p_2p_3=0
  \end{aligned}
  \right.
  \end{equation*}
  and therefore 
  \begin{equation*}
      {(kp_1-p_3)}^2=\rh(k^2p_1+p_3).
  \end{equation*}
  More precisely, Theorem~\ref{mmpttt}\,$(c)$ yields that if $k\ge2$, the $T$ form a parabolic arc and if $k=1$, the
  $T$ form two arcs of a parabola.
\end{proof}

\begin{rem}\label{bin:exp}
  The equality
  \begin{equation}
    \label{extr2}
    \max_x\bigabs{r_1\ei{(t_1+\la_1x)}+r_2\ei{(t_2+\la_2x)}}=r_1+r_2
  \end{equation}
  shows that the exposed points of the unit ball of the space
  $\Cont_{\{\la_1,\la_2\}}$ are the trigonometric monomials
  $\ei\alpha\e_{\la_1}$ and $\ei\alpha\e_{\la_2}$ with $\alpha\in\R$ and that
  no trigonometric binomial is an extreme point of the unit ball of
  $\Cont_\Lambda$.
\end{rem}
\begin{proof}[Proof of Theorem~\ref{extrpts}~$(a)$]
  A linear functional on $\Cont_\Lambda$ extends to a linear
  functional on the space of continuous functions, that is, to a
  measure $\mu$, with same norm.  A measure $\mu$ attains its norm
  only at functions with constant modulus on its support. If $\mu$
  attains its norm on a trigonometric binomial, then it attains its
  norm on a trigonometric monomial because this trigonometric binomial
  is a convex combination of two trigonometric monomials with same
  norm by Equation~\eqref{extr2}. If $\mu$ attains its norm on a
  trigonometric trinomial $T$, there are two cases by Theorem~\ref{mmpttt}\,$(b,c)$:
  \begin{itemize}
  \item $T$ attains its maximum modulus at a unique point modulo
    $2\pi/d$: then the support of $\mu$ has only one point modulo
    $2\pi/d$, so that $\mu$ is the multiple of a Dirac measure and
    attains its norm on any trigonometric monomial.
  \item $T$ attains its maximum modulus at two points modulo
    $2\pi/d$.
  \end{itemize}

  Conversely, the trigonometric monomial $\ei\alpha\e_\la$ is
  exposed to the linear form
  \begin{equation*}
    U\mapsto\frac1{2\pi}\int_0^{2\pi}U(x)\emi{(\alpha+\la x)}\mathrm
  dx.
\end{equation*}
  A trigonometric trinomial $T$ that attains its maximum
  modulus, $1$, at two points $x^*_1$ and $x^*_2$ modulo $2\pi/d$ is
  exposed, by Lemma~\ref{deter}\,$(a)$, to any nontrivial convex
  combination of the unimodular multiples of Dirac measures
  \begin{math}
    \overline{T(x^*_1)}\delta_{x^*_1}
  \end{math}
  and
  \begin{math}
    \overline{T(x^*_2)}\delta_{x^*_2}.
  \end{math}
\end{proof}

\begin{rem}
  This is a complex version of \cite[Lemma~2.3]{re95}.
\end{rem}

\begin{proof}[Proof of Theorem~\ref{extrpts}~$(b)$]
  Let $K$ be the unit ball of $\Cont_\Lambda$. Straszewicz's Theorem
  \cite{st35} tells that the exposed points of $K$ are dense in the set of its
  extreme points. Let $U$ be a limit point of exposed points of $K$. If $U$ is
  a trigonometric monomial, $U$ is exposed. If $U$ is a trigonometric
  binomial, $U$ is not an extreme point of $K$ by Remark~\ref{bin:exp}. If $U$
  is a trigonometric trinomial, it is the limit point of trigonometric
  trinomials that attain their maximum modulus twice modulo $2\pi/d$, so that
  either $U$ also attains its maximum modulus twice modulo $2\pi/d$ or, by
  Rolle's Theorem, $U$ attains its maximum modulus with multiplicity $4$. Let
  us prove that if a trigonometric trinomial $U$ attains its maximum modulus
  with multiplicity $4$ at a point $x$, then $U$ is an extreme point of $K$.
  Suppose that $U$ is the midpoint of two points $A$ and $B$ in $K$. Then
  \begin{math}
    \abs{A(x)}\le1,
  \end{math}
  $\abs{B(x)}\le1$ and $\bigl(A(x)+B(x)\bigr)\big/2=U(x)$, so that
  $A(x)=B(x)=U(x)$. Furthermore
  \begin{align*}
    \abs{U(x+h)}
    &\le\frac{\abs{A(x+h)}+\abs{B(x+h)}}2\\
    &=1+\frac{h^2}4\left(\frac{\mathrm d^2{\abs A}}{\mathrm dx^2}(x)+\frac{\mathrm d^2{\abs B}}{\mathrm dx^2}(x)\right)+o(h^2)
  \end{align*}
  so that, as $\abs{U(x+h)}=1+o(h^3)$,
  \begin{equation*}
    \frac{\mathrm d^2{\abs A}}{\mathrm dx^2}(x)\le0,~\frac{\mathrm d^2{\abs B}}{\mathrm dx^2}(x)\le0~\text{and}~\frac{\mathrm d^2{\abs A}}{\mathrm dx^2}(x)+\frac{\mathrm d^2{\abs B}}{\mathrm dx^2}(x)\ge0
  \end{equation*}
  and therefore $A$ and $B$ also attain their maximum modulus with
  multiplicity $4$ at $x$. As this implies that $A$ and $B$ are trigonometric
  trinomials, Lemma~\ref{deter}\,$(b)$ yields that $U$, $A$ and $B$ lie on a
  parabola: this implies $A=B=U$.
\end{proof}

\begin{rem}
  The set of extreme points of the unit ball of $\Cont_\Lambda$ is not
  closed: for example, if $\la_2$ is between $\la_1$ and $\la_3$,
  every absolutely convex combination of $\e_{\la_1}$ and $\e_{\la_3}$
  is a limit point of exposed points.
\end{rem}

\begin{rem}
  If $\la_1-\la_2$ is not a multiple of $\la_3-\la_2$, nor vice
  versa, then we obtain that every extreme point of the unit ball of
  $\Cont_\Lambda$ is exposed.
\end{rem}

\begin{rem}
  In particular, compare our description of the extreme points of the
  unit ball of $\Cont_{\{0,1,2\}}$ with the characterisation given in \cite{dy03}.
\end{rem}

\section{Dependence of the maximum modulus on the arguments}
\label{seccheb}

We wish to study how the maximum modulus of a trigonometric trinomial
depends on the phase of its coefficients. We shall use the following
formula that gives an expression for the directional derivative of a
maximum function. It was established in \cite{ch43}. Elementary properties of
maximum functions are addressed in \cite[Part Two, Problems 223--226]{ps72}.

\begin{Ch}
  [\protect{\cite[Th.~VI.3.2, (3.6)]{dm90}}]  
  Let $I\subset\R$ be an open interval and let $K$ be a compact space. Let
  $f(t,x)$ be a function on $I\times K$ that is continuous together
  with $\dfrac{\partial f}{\partial t}(t,x)$. Let
  \begin{equation*}
    f^*(t)=\max_{x\in K}f(t,x).
  \end{equation*}
  Then $f^*(t)$ admits the following expansion at every $t\in I$:
  \begin{equation}
    \label{4}
    f^*(t+h)
    =f^*(t)
    +\max_{f(t,x)=f^*(t)}\left(h\dfrac{\partial f}{\partial t}(t,x)\right)
    +o(h).
  \end{equation}
\end{Ch}

\begin{prp}
  Let $k,l$ be two positive coprime integers. Let $r_1$, $r_2$ and
  $r_3$ be three positive real numbers. Then
  \begin{equation*}
    \max_{x}\bigabs{r_1\emi{kx}+r_2\ei{t}+r_3\ei{lx}}
  \end{equation*}
  is an even $2\pi/(k+l)$-periodic function of $t$ that strictly decreases on
  $[0,\pi/(k+l)]$: in particular
  \begin{equation*}
    \min_{t}\max_{x}\bigabs{r_1\emi{kx}+r_2\ei{t}+r_3\ei{lx}}
    =\max_{x}\bigabs{r_1\emi{kx}+r_2\ei{\pi/(k+l)}+r_3\ei{lx}}.
  \end{equation*}
\end{prp}

\begin{proof}
  Let
  \begin{equation}
    \label{f7}
    f(t,x)={\bigabs{r_1\emi{kx}+r_2\ei{t}+r_3\ei{lx}}}^2.
  \end{equation}
  By~\eqref{pair} and~\eqref{per}, $f^*$ is an even
  $2\pi/(k+l)$-periodic function. 

  Let $t\in\io[0,\pi/(k+l)]$ and choose $x^*$ such that $f(t,x^*)=f^*(t)$: then $x^*\in[-t/k,t/l]$ by
  Proposition~\ref{loc}, so that
  \begin{equation*}
    \frac1{2r_2}\dfrac{\partial f}{\partial t}(t,x^*)
    =-r_1\sin(t+kx^*)-r_3\sin(t-lx^*)<0
  \end{equation*}
  because $t+kx^*\in[0,(k+l)t/l]$ and $t-lx^*\in[0,(k+l)t/k]$ do
  not vanish simultaneously. By Formula~\eqref{4}, $f^*$ decreases
  strictly on $[0,\pi/(k+l)]$.
\end{proof}

\begin{prp}
  \label{croissant}
  Let $k,l$ be two positive coprime integers. Let $r_1$, $r_2$ and
  $r_3$ be three positive real numbers. Then
  \begin{equation}
    \label{gg}
    \frac
    {\max_{x}\bigabs{r_1\emi{kx}+r_2\ei{t}+r_3\ei{lx}}}
    {\bigabs{r_1+r_2\ei{t}+r_3}}
  \end{equation}
  is an increasing function of $t\in[0,\pi/(k+l)]$. If $kr_1=lr_3$, it is
  constantly equal to $1$; otherwise it is strictly increasing. 
\end{prp}

\begin{proof}
  Let $f(t,x)$ be as
  in~\eqref{f7}: then the expression~\eqref{gg} is $g^*(t)^{1/2}$ with
  \begin{gather*}
    g(t,x)=\frac{f(t,x)}{f(t,0)}.
  \end{gather*}
  If $kr_1=lr_3$, then $f(t,0)=f^*(t)$, so that $g^*(t)=1$. As shown in
  the beginning of Section~\ref{sec4}, we may suppose without loss of
  generality that $kr_1<lr_3$.  Let $t\in\io[0,\pi/(k+l)]$ and choose $x^*$
  such that $f(t,x^*)=f^*(t)$: then $x^*\in\io[0,t/l]$ by
  Propositions~\ref{loc} and~\ref{prp-2} and
  \begin{multline*}
    \frac{f(t,0)^2}{2r_2}\dfrac{\partial g}{\partial t}(t,x^*)
    =\frac1{2r_2}\biggl(
    \dfrac{\partial f}{\partial t}(t,x^*)f(t,0)
    -f(t,x^*)\dfrac{\partial f}{\partial t}(t,0)\biggr)\\
    \begin{aligned}
      &=\bigl(-r_1\sin(t+kx^*)-r_3\sin(t-lx^*)\bigr)f(t,0)+
      f^*(t)(r_1+r_3)\sin t\\
      &=h(0)f^*(t)-h(x^*)f(t,0)
    \end{aligned}
  \end{multline*}
  with 
  \begin{equation*}
    h(x)=r_1\sin(t+kx)+r_3\sin(t-lx).
  \end{equation*}
  Let us show that $h$ is strictly decreasing on
  $[0,t/l]$: in fact, if $x\in\io[0,t/l]$,
  \begin{equation*}
    \frac{\mathrm dh}{\mathrm dx}(x)=kr_1\cos(t+kx)-lr_3\cos(t-lx)
    <(kr_1-lr_3)\cos(t-lx)<0.
  \end{equation*}
  As $f^*(t)>f(t,0)$ and $h(0)>h(x^*)$, $({\partial g}/{\partial
    t})(t,x^*)>0$. By N.~G. Che\-bo\-ta\-r\"ev's formula, $g^*$ increases
  strictly on $[0,\pi/(k+l)]$.
\end{proof}

It is possible to describe the decrease of the maximum modulus of a
trigonometric trinomial independently of the $r$'s as follows.
\begin{prp}
  \label{prp-5}
  Let $k$ and $l$ be two positive coprime integers. Let $r_1$, $r_2$
  and $r_3$ be three positive real numbers. Let $0\le t'<t\le\pi/(k+l)$.
  Then
  \begin{equation}
    \label{ineg-4}
    \max_x\bigabs{r_1\emi{kx}+r_2\ei{t'}+r_3\ei{lx}}
    \le\frac{\cos(t'/2)}{\cos(t/2)}
    \max_x\bigabs{r_1\emi{kx}+r_2\ei{t}+r_3\ei{lx}}
  \end{equation}
  with equality if and only if $r_1:r_2:r_3=l:k+l:k$.
\end{prp}
\begin{proof}
  Let us apply Proposition~\ref{croissant}. We have
  \begin{align*}
    \frac{{\bigabs{r_1+r_2\ei{t'}+r_3}}^2} {{\bigabs{r_1+r_2\ei{t}+r_3}}^2}
    &=1+\frac{2r_2(r_1+r_3)(\cos t'-\cos t)}{(r_1+r_3)^2+2r_2(r_1+r_3)\cos t+r_2^2}\\
    &=1+\frac{\cos t'-\cos t}{\cos t+\bigl(r_2^2+(r_1+r_3)^2\bigr)/2r_2(r_1+r_3)}\\
    &\le1+\frac{\cos t'-\cos t}{\cos t+1}=\frac{\cos t'+1}{\cos t+1}
  \end{align*}
  by the arithmetic-geometric inequality, with equality if and only if $r_2=r_1+r_3$.
  Therefore Inequality~\eqref{ineg-4} holds, with equality if and only if
  $kr_1=lr_3$ and $r_2=r_1+r_3$.
\end{proof}

We may now find the minimum of the maximum modulus of a trigonometric
trinomial with fixed Fourier coefficient arguments and moduli sum.
Proposition~\ref{prp-5} yields with $t'=0$

\begin{cor}
  \label{prp-3}
  Let $k$ and $l$ be two positive coprime integers. Let $r_1$, $r_2$
  and $r_3$ be three positive real numbers. Let $t\in\iog[0,\pi/(k+l)]$.
  Then
  \begin{equation*}
    \frac{\max_x\bigabs{r_1\emi{kx}+r_2\ei{t}+r_3\ei{lx}}}{r_1+r_2+r_3}\ge\cos(t/2)
  \end{equation*}
  with equality if and only if $r_1:r_2:r_3=l:k+l:k$.
\end{cor}

\begin{rem}
  There is a shortcut proof of Corollary~\ref{prp-3}:
  \begin{align*}
    \frac{\max_x\bigabs{r_1\emi{kx}+r_2\ei{t}+r_3\ei{lx}}}{r_1+r_2+r_3}
    &\ge\frac{\bigabs{r_1+r_2\ei{t}+r_3}}{r_1+r_2+r_3}\\
    &=\sqrt{1-\frac{4(r_1+r_3)r_2}{{(r_1+r_2+r_3)}^2}\sin^2(t/2)}\\
    &\ge\sqrt{1-\sin^2(t/2)}=\cos(t/2)
  \end{align*}
  and equality holds if and only if
  $\bigabs{r_1\emi{kx}+r_2\ei{t}+r_3\ei{lx}}$ is maximal for $x=0$ and
  $r_1+r_3=r_2$.
\end{rem}

\section{\sloppy The norm of unimodular relative Fourier multipliers}
\label{secnm}

We may now compute the norm of unimodular relative Fourier multipliers. 

\begin{cor}
  Let $k$ and $l$ be two positive coprime integers. Let
  $t\in[0,\pi/\allowbreak(k+l)]$. Let $M$ be the relative Fourier multiplier
  $(0,t,0)$   that maps the element
  \begin{equation}
    \label{Tkl}
    r_1\ei{u_1}\e_{-k}+r_2\ei{u_2}\e_0+r_3\ei{u_3}\e_{l}
  \end{equation}
  of the normed space $\Cont_{\{-k,0,l\}}$ on
  \begin{equation*}
    r_1\ei{u_1}\e_{-k}+r_2\ei{(t+u_2)}\e_0+r_3\ei{u_3}\e_{l}.
  \end{equation*}
  Then $M$ has norm
  $\cos\bigl(\pi/2(k+l)-t/2\bigr)\big/\cos\bigl(\pi/2(k+l)\bigr)$
  and attains its norm exactly at elements of form~\eqref{Tkl}
  with $r_1:r_2:r_3=l:k+l:k$ and
  \begin{equation*}
    -lu_1+(k+l)u_2-ku_3=\pi\mod{2\pi}.
  \end{equation*}
\end{cor}

\begin{proof}
  This follows from Proposition~\ref{prp-5} and the concavity of $\cos$ on $[0,\pi/2]$.
\end{proof}

\begin{rem}
  This corollary permits to guess how to lift $M$ to an operator that acts by
  convolution with a measure $\mu$. Note that $\mu$ is a Hahn-Banach extension
  of the linear form $f\mapsto Mf(0)$. The relative multiplier $M$ is an isometry if
  and only if $t=0$ and $\mu$ is the Dirac measure in $0$. Otherwise, $t\ne0$;
  the proof of Theorem~\ref{extrpts}\,$(a)$ shows that $\mu$ is a linear
  combination $\alpha\delta_y+\beta\delta_w$ of two Dirac measures such that the
  norm of $M$ is $\abs{\alpha}+\abs{\beta}$. Let
  $f(x)=l\emi{kx}+(k+l)\ei{\pi/(k+l)}+k\ei{lx}$: $M$ attains its norm at $f$,
  $f$ attains its maximum modulus at $0$ and $2m\pi/(k+l)$, and $Mf$ attains
  its maximum modulus at $2m\pi/(k+l)$, where $m$ is the inverse of $l$ modulo
  $k + l$. As
  \begin{align*}
    (\abs{\alpha}+\abs{\beta})\max_x\abs{f(x)} 
    &=\max_x\abs{Mf(x)}\\
    &=\abs{\mu*f\bigl(2m\pi/(k+l)\bigr)}\\
    &=\abs{\alpha f\bigl(2m\pi/(k+l)-y\bigr)+\beta
      f\bigl(2m\pi/(k+l)-w\bigr)},
  \end{align*}
  we must choose $\{y,w\}=\{0,2m\pi/(k+l)\}$. A computation yields then
  \begin{equation*}
    \mu=\ei{t/2}\frac{\sin\bigl(\pi/(k+l)-t/2\bigr)}{\sin\bigl(\pi/(k+l)\bigr)}\delta_0
    +\ei{(t/2+\pi/(k+l))}\frac{\sin(t/2)}{\sin\bigl(\pi/(k+l)\bigr)}\delta_{2m\pi/(k+l)}.
  \end{equation*}
  Consult \cite{sh61} on this issue. 
\end{rem}

\section{The Sidon constant of integer sets}
\label{secsidon}

Let us study the maximum modulus of a trigonometric trinomial with given
Fourier coefficient moduli sum. We get the following result as an immediate
consequence of Corollary~\ref{prp-3}. 
\begin{prp}
  \label{prp-4}
  Let $k$ and $l$ be two positive coprime integers. Let $r_1$, $r_2$
  and $r_3$ be three positive real numbers. Let $t\in[0,\pi/(k+l)]$.
  Then
  \begin{equation*}
    \max_x\bigabs{r_1\emi{kx}+r_2\ei{t}+r_3\ei{lx}}
    \ge(r_1+r_2+r_3)\cos\bigl(\pi/2(k+l)\bigr)
  \end{equation*}
  with equality if and only if $r_1:r_2:r_3=l:k+l:k$ and $t=\pi/(k+l)$.
\end{prp}

This means that the Sidon constant of $\{-k,0,l\}$ equals
$\sec\bigl(\pi/2(k+l)\bigr)$.

The Sidon constant of integer sets was known in three instances
only:
\begin{itemize}
\item The equality
  \begin{equation*}
    \max_x\bigabs{r_1\ei{(t_1+\la_1x)}+r_2\ei{(t_2+\la_2x)}}=r_1+r_2
  \end{equation*}
  shows that the Sidon constant of sets with one or two elements is $1$.
\item The Sidon constant of $\{-1,0,1\}$ is $\sqrt2$ and it is attained for
  $\e_{-1}+2\iu+\e_{1}$. Let us give the original argument: if
  $f(x)=\bigabs{r_1\emi{x}+r_2\ei{t}+r_3\ei{x}}^2$, the parallelogram
  identity and the arithmetic-quadratic inequality yield
  \begin{multline*}
    \max_x{f(x)}
    \ge\max_x\frac{f(x)+f(x+\pi)}2\\
    \begin{aligned}
      &=\max_x\frac{
        {\bigabs{r_1\emi{x}+r_3\ei{x}+r_2\ei{t}}}^2+
        {\bigabs{r_1\emi{x}+r_3\ei{x}-r_2\ei{t}}}^2}2\\
      &=\max_x{\bigabs{r_1\emi{x}+r_3\ei{x}}}^2+{\bigabs{r_2\ei{t}}}^2\\
      &=(r_1+r_3)^2+r_2^2\ge\frac{(r_1+r_2+r_3)^2}2.
    \end{aligned}
  \end{multline*}
\item The Sidon constant of $\{0,1,2,3,4\}$ is $2$ and it is attained for
  $1+2\e_1+2\e_2-2\e_3+\e_4$.
\end{itemize}
These results were obtained by D.~J. Newman (see \cite{sh51}.)  The fact that
the Sidon constant of integer sets with three elements cannot be $1$ had been
noted with pairwise different proofs in \cite{sh51,chm81,li96}.\medskip

The following estimates for the Sidon constant of large integer sets are known.
\begin{itemize}
\item E.~Beller and D.~J. Newman \cite{bn71} showed that the Sidon constant of
  $\{0,\allowbreak 1,\allowbreak \dots,\allowbreak n\}$ is equivalent to $\sqrt n$.
\item (Hadamard sets.) Let $q>1$ and suppose that the sequence $(\lambda_j)_{j\ge1}$ grows with
  geometric ratio $q$: $\abs{\la_{j+1}}\ge q\abs{\la_j}$ for every $j$. Then
  the Sidon constant of $\{\la_1,\la_2,\dots\}$ is finite; it is at most $8$
  if $q\ge2$, it is at most $2$ if $q\ge3$ (see \cite{lr75}), and it is at
  most $1+\pi^2\big/\left(2q^2-2-\pi^2\right)$ if $q>\sqrt{1+\pi^2/2}$ (see
  \cite[Corollary 9.4]{ne98} or the updated \cite[Corollary 10.2.1]{ne99}.)
\end{itemize}

Our computations show that the last estimate of the Sidon constant has the
right order in $q^{-1}$ for geometric progressions.

\begin{prp}
  Let $C$ be the Sidon constant of the geometric progression
  $\{1,q,q^2,\dots\}$, where $q\ge3$ is an integer. Then
  \begin{equation*}
    1+\pi^2\big/8(q+1)^2
    \le\sec\bigl(\pi/2(q+1)\bigr)
    \le C\le1+\pi^2\big/\left(2q^2-2-\pi^2\right).
  \end{equation*}
\end{prp}

One initial motivation for this work was to decide whether there are
sets $\{\lambda_j\}_{j\ge1}$ with $\abs{\la_{j+1}}\ge q\abs{\la_j}$
whose Sidon constant is arbitrarily close to $1$ and to find evidence
among sets with three elements. That there are such sets,
arbitrarily large albeit finite, may in fact be proved by the method of
Riesz products in \cite[Appendix V, \S1.II]{ks63}; see also
\cite[Proposition 13.1.3]{ne99}. The case of infinite sets
remains open.

A second motivation was to show that the real and complex unconditional
constants of the basis $(\e_{\la_1},\e_{\la_2},\e_{\la_3})$ of
$\Cont_\Lambda$ are different; we prove however that they coincide, and it remains an open
question whether they may be different for larger sets.  The \emph{real
  unconditional constant} of $(\e_{\la_1},\e_{\la_2},\e_{\la_3})$ is the
maximum of the norm of the eight unimodular relative Fourier multipliers
$(t_1,t_2,t_3)$ such that $t_k=0$ modulo $\pi$. Let $i,j,k$ be a permutation
of $1,2,3$ such that the power of $2$ in $\la_i-\la_k$ and in $\la_j-\la_k$
are equal. Lemma~\ref{isom} shows that the four relative multipliers satisfying
$t_i=t_j$ modulo $2\pi$ are isometries and that the norm of any of the four
others, satisfying $t_i\ne t_j$ modulo $2\pi$, gives the real unconditional
constant. In general, the complex unconditional constant is bounded by $\pi/2$
times the real unconditional constant, as proved in \cite{se97}; in our case,
they are equal.
\begin{cor}
  The complex unconditional constant of the basis
  $(\e_{\la_1},\allowbreak\e_{\la_2},\allowbreak\e_{\la_3})$ of $\Cont_\Lambda$ is equal to its real
  unconditional constant.
\end{cor}

\begin{acknowledgements}
  The Hungarian-French Scientific and Tech\-no\-lo\-gi\-cal Go\-vern\-men\-tal Cooperation,
  Project \#F-10/04, eased scientific contacts that were helpful to this work.
\end{acknowledgements}

\def\cprime{$'$}

\begin{sloppypar}
  \noindent \textit{Keywords.} Trigonometric trinomial, maximum
    modulus, exposed point, ex\-treme point, Mandel$'$shtam problem, extremal
    problem, relative Fourier multiplier, Sidon constant, unconditional
    constant.  
\end{sloppypar}
\medskip

\noindent 
2000 \textit{Mathematics Subject Classification.} Primary 30C10, 42A05,
42A45, 46B20; Secondary 26D05, 42A55, 46B15.
\medskip

\noindent
Stefan Neuwirth, Laboratoire de Math\'ematiques, Universit\'e de
Franche-Comt\'e, 25030 Besan\c con cedex, France,
{stefan.neuwirth@univ-fcomte.fr}.


\begin{thebibliography}{10}

\bibitem{bn71}
E.~Beller and D.~J. Newman, {\em An {$l\sb{1}$} extremal problem for
  polynomials}, Proc. Amer. Math. Soc. 29 (1971),  474--481.

\bibitem{chm81}
D.~I. Cartwright, R.~B. Howlett and J.~R. McMullen, {\em Extreme values for the
  {S}idon constant}, Proc. Amer. Math. Soc. 81 (1981),  531--537.

\bibitem{ch43}
N.~G. Chebotar\"ev, {\em On a general criterion of the minimax}, C. R.
  (Doklady) Acad. Sci. URSS (N.S.) 39 (1943),  339--341.

\bibitem{ce49}
\leavevmode\vrule height 2pt depth -1.6pt width 23pt, {\em Kriteri{\u\i}
  minimaksa i ego prilozheniya}, in: Sobranie sochineni\u\i. {V}ol. 2,
  Izdatel\cprime stvo Akademii Nauk SSSR, Moscow-Leningrad, 1949,  396--409.

\bibitem{dm90}
V.~F. Dem{\cprime}yanov and V.~N. Maloz{\"e}mov, {\em Introduction to minimax},
  Halsted Press [John Wiley \& Sons], New York-Toronto, 1974.
\newblock Translated from the Russian by D. Louvish.

\bibitem{dy03}
K.~M. Dyakonov, {\em Extreme points in spaces of polynomials}, Math. Res. Lett.
  10 (2003),  717--728.

\bibitem{ha87}
S.~Hartman, {\em Some problems and remarks on relative multipliers}, Colloq.
  Math. 54 (1987),  103--111.
\newblock Corrected in Colloq. Math. 57 (1989), 189.

\bibitem{ks63}
J.-P. Kahane and R.~Salem, {\em Ensembles parfaits et s\'eries
  tri\-go\-no\-m\'e\-triques}, Ac\-tua\-li\-t\'es Sci. Indust. 1301, Hermann,
  Paris, 1963.

\bibitem{li96}
D.~Li, {\em Complex unconditional metric approximation property for {${\mathscr
  C}\sb \Lambda({\mathbb T})$} spaces}, Studia Math. 121 (1996),  231--247.

\bibitem{li24}
J.~E. Littlewood, {\em A theorem on power series}, {Proc. London Math. Soc.
  (2)} 23 (1925),  94--103.

\bibitem{lr75}
J.~M. L{\'o}pez and K.~A. Ross, {\em Sidon sets}, Lecture Notes Pure Appl.
  Math. 13, Marcel Dekker Inc., New York, 1975.

\bibitem{ne01b}
S.~Neuwirth, {\em The {S}idon constant of sets with three elements}.
\newblock \texttt{arxiv.org/\allowbreak math/\allowbreak 0102145}, 2001.

\bibitem{ne98}
\leavevmode\vrule height 2pt depth -1.6pt width 23pt, {\em Metric
  unconditionality and {F}ourier analysis}, Studia Math. 131 (1998),  19--62.

\bibitem{ne99}
\leavevmode\vrule height 2pt depth -1.6pt width 23pt, {\em Multiplicateurs et
  analyse fonctionnelle}, PhD thesis, Universit\'e Paris 6, 1999.
\newblock \texttt{tel.archives-ouvertes.fr/tel-00010399}.

\bibitem{ps72}
G.~P\'olya and G.~Szeg{\H{o}}, {\em Problems and theorems in analysis. {V}ol.
  {I}: {S}eries, integral calculus, theory of functions}, Grundlehren math.
  Wiss. 193, Springer-Verlag, New York, 1972.
\newblock Translated from the German by D. Aeppli.

\bibitem{re95}
Sz.~Gy. R{\'e}v{\'e}sz, {\em Minimization of maxima of nonnegative and positive
  definite cosine polynomials with prescribed first coefficients}, Acta Sci.
  Math. (Szeged) 60 (1995),  589--608.

\bibitem{sa33}
R.~Salem, {\em Sur les propri{\'e}t{\'e}s extr{\'e}males de certains polynomes
  tri\-go\-no\-m{\'e}\-triques}, C. R. Acad. Sci. Paris 196 (1933),
  1776--1778.

\bibitem{se97}
J.~A. Seigner, {\em Rademacher variables in connection with complex scalars},
  Acta Math. Univ. Comenian. (N.S.) 66 (1997),  329--336.

\bibitem{sh51}
H.~S. Shapiro, {\em Extremal problems for polynomials and power series},
  Master's thesis, Massachusetts Institute of Technology, 1951.

\bibitem{sh61}
\leavevmode\vrule height 2pt depth -1.6pt width 23pt, {\em On a class of
  extremal problems for polynomials in the unit circle}, Portugal. Math. 20
  (1961),  67--93.

\bibitem{st35}
S.~Straszewicz, {\em {\"U}ber exponierte {P}unkte abgeschlossener
  {P}unktmengen}, Fund. Math. 24 (1935),  139--143.

\end{thebibliography}
\end{document}